\documentclass[english]{amsart}
\usepackage[T1]{fontenc}
\usepackage[latin1]{inputenc}
\usepackage{amssymb}

\makeatletter
 \theoremstyle{plain}    
 
 \numberwithin{equation}{section} 
 \numberwithin{figure}{section} 
 \theoremstyle{plain}
 \theoremstyle{plain}    
 \newtheorem*{thm*}{Theorem} 
 \theoremstyle{remark}    
 \newtheorem*{conclusion*}{Conclusion} 

\usepackage{babel}
\makeatother
\begin{document}

\title{On A Relation Between $q$-Exponential And $\theta$-Function}

\author{Ruiming Zhang}

\date{June 6, 2006}

\address{School of Mathematics\\
Guangxi Normal University\\
Guilin City, Guangxi 541004\\
P. R. China.}

\email{ruimingzhang@yahoo.com}



\begin{abstract}
We will use a discrete analogue of the classical Laplace method to
show that for infinitely many positive integers $n$, the main term
of the asymptotic expansion of the scaled $q$-exponential $(-q^{-nt+1/2}u;q)_{\infty}$
could be expressed in $\theta$-function. 
\end{abstract}
\maketitle

\section{Introduction}

Given an arbitrary positive real number $0<q<1$ and an arbitrary
complex number $a$, we define\cite{Gasper,Ismail2} \begin{equation}
(a;q)_{\infty}=\prod_{k=0}^{\infty}(1-aq^{k})\label{eq:1.1}\end{equation}
 and the $q$-shifted factorial as\begin{equation}
(a;q)_{n}=\frac{(a;q)_{\infty}}{(aq^{n};q)_{\infty}}\label{eq:1.2}\end{equation}
 for any integer $n$. Assume that $|z|<1$, the $q$-Binomial theorem
is\cite{Gasper,Ismail2} \begin{equation}
\frac{(az;q)_{\infty}}{(z;q)_{\infty}}=\sum_{k=0}^{\infty}\frac{(a;q)_{k}}{(q;q)_{k}}z^{k}.\label{eq:1.3}\end{equation}
 This could be seen easily from the fact that \begin{equation}
\frac{(az;q)_{\infty}}{(z;q)_{\infty}}=\prod_{k=0}^{\infty}\frac{1-azq^{k}}{1-zq^{k}}\label{eq:1.4}\end{equation}
 defines an analytic function in the region $|z|<1$. Hence it could
be expanded as a power series \begin{equation}
\frac{(az;q)_{\infty}}{(z;q)_{\infty}}=\sum_{k=0}^{\infty}a_{k}z^{k}\label{eq:1.5}\end{equation}
 for $|z|<1$ with \begin{equation}
a_{0}=1.\label{eq:1.6}\end{equation}
Observe that \begin{equation}
\frac{(az;q)_{\infty}}{(z;q)_{\infty}}=\frac{1-az}{1-z}\frac{(azq;q)_{\infty}}{(zq;q)_{\infty}},\label{eq:1.7}\end{equation}
 thus we have\begin{equation}
(1-z)\sum_{k=0}^{\infty}a_{k}z^{k}=(1-az)\sum_{k=0}^{\infty}a_{k}q^{k}z^{k}\label{eq:1.8}\end{equation}
 or\begin{equation}
\sum_{k=0}^{\infty}a_{k}(1-q^{k})z^{k}=\sum_{k=1}^{\infty}a_{k-1}(1-aq^{k-1})z^{k}.\label{eq:1.9}\end{equation}
 Then \begin{equation}
a_{k}=\frac{1-aq^{k-1}}{1-q^{k}}a_{k-1}\label{eq:1.10}\end{equation}
 and\begin{equation}
a_{k}=\frac{(a;q)_{k}}{(q;q)_{k}}\label{eq:1.11}\end{equation}
 for $k=1,2,\dotsc$. For any complex number $z$, let $a$ be a large
real number with $a\gg|z|$. We replace $z$ in (\ref{eq:1.3}) by
$\frac{z}{a}$, \begin{equation}
\frac{(z;q)_{\infty}}{(z/a;q)_{\infty}}=\sum_{k=0}^{\infty}\frac{(a;q)_{k}/a^{k}}{(q;q)_{k}}z^{k}.\label{eq:1.12}\end{equation}
 Since \begin{equation}
\lim_{a\to\infty}\frac{(a;q)_{k}}{a^{k}}=(-1)^{k}q^{k(k-1)/2},\label{eq:1.13}\end{equation}
 then\begin{equation}
\lim_{a\to\infty}\sum_{k=0}^{\infty}\frac{(a;q)_{k}/a^{k}}{(q;q)_{k}}z^{k}=\sum_{k=0}^{\infty}\frac{q^{k(k-1)/2}}{(q;q)_{k}}(-z)^{k}\label{eq:1.14}\end{equation}
 by the Lebesgue dominated convergent theorem. Similarly from (\ref{eq:1.3})\begin{equation}
\frac{1}{(z/a;q)_{\infty}}=\sum_{k=0}^{\infty}\frac{1}{(q;q)_{k}}\left(\frac{z}{a}\right)^{k},\label{eq:1.15}\end{equation}
 then\begin{eqnarray}
\lim_{a\to\infty}\frac{(z;q)_{\infty}}{(z/a;q)_{\infty}} & = & \lim_{a\to\infty}\sum_{k=0}^{\infty}\frac{(z;q)_{\infty}}{(q;q)_{k}}\left(\frac{z}{a}\right)^{k}\label{eq:1.16}\\
 & = & (z;q)_{\infty}\nonumber \end{eqnarray}
 by the Lebesgue dominated theorem. Thus we have proved that for any
complex number $z$\begin{equation}
(z;q)_{\infty}=\sum_{k=0}^{\infty}\frac{q^{k(k-1)/2}}{(q;q)_{k}}(-z)^{k}.\label{eq:1.17}\end{equation}
 Notice that\begin{equation}
((1-q)z;q)_{\infty}=\sum_{k=0}^{\infty}\frac{(1-q)^{k}}{(q;q)_{k}}q^{k(k-1)/2}(-z)^{k}.\label{eq:1.18}\end{equation}
 Since\begin{equation}
\frac{1-q^{k}}{1-q}\ge kq^{k-1}\label{eq:1.19}\end{equation}
 for $k\ge1$, thus\begin{equation}
\frac{(1-q)^{k}q^{k(k-1)/2}}{(q;q)_{k}}\le\frac{1}{k!}\label{eq:1.20}\end{equation}
for $k\ge1$, hence by the Lebesgue dominated theorem we have\begin{equation}
\lim_{q\to1}((1-q)z;q)_{\infty}=\sum_{k=0}^{\infty}\frac{(-z)^{k}}{k!}=e^{-z},\label{eq:1.21}\end{equation}
and this is the reason why $(z;q)_{\infty}$ is called a $q$-exponential.
Indeed, it is one of several $q$-analogues of $e^{z}$ in the theory
of $q$-series. From (\ref{eq:1.18}) and (\ref{eq:1.20}), it is
also clear that \begin{equation}
|((1-q)z;q)_{\infty}|\le e^{|z|}\label{eq:1.22}\end{equation}
for any complex number $z$. For any nonzero complex number $z$,
we define the theta function as \begin{equation}
\theta(z;q)=\sum_{k=-\infty}^{\infty}q^{k^{2}/2}z^{k}.\label{eq:1.23}\end{equation}
 The Jacobi's triple product formula says that\cite{Gasper,Ismail2}
\begin{equation}
\sum_{k=-\infty}^{\infty}q^{k^{2}/2}z^{k}=(q,-q^{1/2}z,-q^{1/2}/z;q)_{\infty}.\label{eq:1.24}\end{equation}
 For any positive real number $t$, we consider the following set
\begin{equation}
\mathbb{S}(t)=\left\{ \{ nt\}:n\in\mathbb{N}\right\} .\label{eq:1.25}\end{equation}
 It is clear that $\mathbb{S}(t)\subset[0,1)$ and it is a finite
set when $t$ is a positive rational number. In this case, for any
$\lambda\in\mathbb{S}(t)$, there are infinitely many positive integers
$n$ and $m$ such that \begin{equation}
nt=m+\lambda,\label{eq:1.26}\end{equation}
where\begin{equation}
m=\left\lfloor nt\right\rfloor .\label{eq:1.27}\end{equation}
If $t$ is a positive irrational number, then $\mathbb{S}(t)$ is
a subset of $(0,1)$ with infinite elements, and it is well-known
that $\mathbb{S}(t)$ is uniformly distributed in $(0,1)$. A theorem
of Chebyshev\cite{Hua} says that given any $\beta\in[0,1)$, there
are infinitely many positive integers $n$ and $m$ such that\begin{equation}
nt=m+\beta+\gamma_{n}\label{eq:1.28}\end{equation}
 with\begin{equation}
|\gamma_{n}|\le\frac{3}{n},\label{eq:1.29}\end{equation}
 It is clear that for $n$ large enough, we may also have\begin{equation}
m=\left\lfloor nt\right\rfloor .\label{eq:1.30}\end{equation}

\section{Main Results}

From (\ref{eq:1.17}), clearly we have\begin{eqnarray}
(-q^{-nt+1/2}u;q)_{\infty} & = & \sum_{k=0}^{\infty}\frac{q^{k^{2}/2-knt}}{(q;q)_{k}}u^{k}.\label{eq:2.1}\end{eqnarray}

\begin{thm*}
Given an arbitrary positive real number $0<q<1$ and an arbitrary
nonzero complex number $u$, we have
\begin{enumerate}
\item For any positive rational number $t$, let $\lambda\in\mathbb{S}(t)$,
there are infinitely many positive integers $n$ and $m$ such that
\begin{equation}
tn=m+\lambda,\label{eq:2.2}\end{equation}
 and for these $n$'s and $m$'s we have\begin{equation}
(-q^{-nt+1/2}u;q)_{\infty}=\frac{u^{m}\left\{ \theta(u^{-1}q^{\lambda};q)+r(n)\right\} }{(q;q)_{\infty}q^{m^{2}/2+m\lambda}}\label{eq:2.3}\end{equation}
 with\begin{eqnarray}
|r(n)| & \le & \frac{3(-q;q)_{\infty}\theta(\left|u\right|^{-1}q^{\lambda};q)}{1-q}\nonumber \\
 & \times & \left\{ q^{m/2}+\frac{q^{m^{2}/8}}{|u|^{\left\lfloor m/2\right\rfloor +1}}\right\} .\label{eq:2.4}\end{eqnarray}
 for $n$, $m$ and $\lambda$ satisfying (\ref{eq:2.2}) with $n$
and $m$ are large enough.
\item Given a positive irrational number $t$, for any real number $\beta\in[0,1)$,
there are infinitely many positive integers $n$ and $m$ such that
\begin{equation}
nt=m+\beta+\gamma_{n}\label{eq:2.5}\end{equation}
 with\begin{equation}
|\gamma_{n}|\le\frac{3}{n}.\label{eq:2.6}\end{equation}
 For such $n$'s and $m$'s we have\begin{equation}
(-q^{-nt+1/2}u;q)_{\infty}=\frac{u^{m}\left\{ \theta(u^{-1}q^{\beta};q)+O\left(\frac{\log n}{n}\right)\right\} }{(q;q)_{\infty}q^{mnt-m^{2}/2}}.\label{eq:2.7}\end{equation}
 for $n$ is large enough.
\end{enumerate}
\end{thm*}
\begin{proof}
In the case that $t$ is a positive rational number, for any $\lambda\in\mathbb{S}(t)$
and $n$, $m$ are large, we have\begin{eqnarray}
(-q^{-nt+1/2}u;q)_{\infty}(q;q)_{\infty} & = & \sum_{k=0}^{\infty}(q^{k+1};q)_{\infty}q^{k^{2}/2-km-k\lambda}u^{k}\label{eq:2.8}\\
 & = & s_{1}+s_{2}\nonumber \end{eqnarray}
 where\begin{equation}
s_{1}=\sum_{k=0}^{m}(q^{k+1};q)_{\infty}q^{k^{2}/2-km-k\lambda}u^{k}\label{eq:2.9}\end{equation}
 and \begin{equation}
s_{2}=\sum_{k=m+1}^{\infty}(q^{k+1};q)_{\infty}q^{k^{2}/2-km-k\lambda}u^{k}.\label{eq:2.10}\end{equation}
 In $s_{1}$ we change the summation from $k$ to $m-k$, we have\begin{eqnarray*}
\frac{s_{1}q^{m^{2}/2+m\lambda}}{u^{m}} & = & \sum_{k=0}^{m}(q^{m-k+1};q)_{\infty}q^{k^{2}/2}(u^{-1}q^{\lambda})^{k}\\
 & = & \sum_{k=0}^{\infty}q^{k^{2}/2}(u^{-1}q^{\lambda})^{k}\\
 & - & \sum_{k=\left\lfloor m/2\right\rfloor +1}^{\infty}q^{k^{2}/2}(u^{-1}q^{\lambda})^{k}\\
 & + & \sum_{k=0}^{\left\lfloor m/2\right\rfloor }q^{k^{2}/2}(u^{-1}q^{\lambda})^{k}\left((q^{m-k+1};q)_{\infty}-1\right)\\
 & + & \sum_{k=\left\lfloor m/2\right\rfloor +1}^{m}(q^{m-k+1};q)_{\infty}q^{k^{2}/2}(u^{-1}q^{\lambda})^{k}\\
 & = & \sum_{k=0}^{\infty}q^{k^{2}/2}(u^{-1}q^{\lambda})^{k}+s_{11}+s_{12}+s_{13}.\end{eqnarray*}
 Clearly,\begin{eqnarray*}
|s_{11}+s_{13}| & \le & 2\sum_{k=\left\lfloor m/2\right\rfloor +1}^{\infty}q^{k^{2}/2}(\left|u\right|^{-1}q^{\lambda})^{k}\\
 & \le & \frac{2q^{m^{2}/8}}{|u|^{\left\lfloor m/2\right\rfloor +1}}\theta(\left|u\right|^{-1}q^{\lambda};q).\end{eqnarray*}
 By (\ref{eq:1.3}), for $0\le k\le\left\lfloor m/2\right\rfloor $,
we have\begin{equation}
\left|(q^{m-k+1};q)_{\infty}-1\right|=\frac{(-q;q)_{\infty}}{1-q}q^{m/2},\label{eq:2.11}\end{equation}
 then\begin{eqnarray*}
|s_{12}| & \le & \frac{(-q;q)_{\infty}}{1-q}q^{m/2}\sum_{k=0}^{\infty}q^{k^{2}/2}(\left|u\right|^{-1}q^{\lambda})^{k}\\
 & \le & \frac{q^{m/2}(-q;q)_{\infty}\theta(\left|u\right|^{-1}q^{\lambda};q)}{1-q},\end{eqnarray*}
 hence\begin{equation}
\frac{s_{1}q^{m^{2}/2+m\lambda}}{u^{m}}=\sum_{k=0}^{\infty}q^{k^{2}/2}(u^{-1}q^{\lambda})^{k}+r_{1}(n)\label{eq:2.12}\end{equation}
 with\begin{eqnarray}
|r_{1}(n)| & \le & \frac{2(-q;q)_{\infty}\theta(\left|u\right|^{-1}q^{\lambda};q)}{1-q}\label{eq:2.13}\\
 & \times & \left\{ q^{m/2}+\frac{q^{m^{2}/8}}{|u|^{\left\lfloor m/2\right\rfloor +1}}\right\} .\nonumber \end{eqnarray}
 In $s_{2}$ we change the summation from $k$ to $k+m$\begin{eqnarray*}
\frac{s_{2}q^{m^{2}/2+m\lambda}}{u^{m}} & = & \sum_{k=1}^{\infty}(q^{m+k+1};q)_{\infty}q^{k^{2}/2}(uq^{-\lambda})^{k}\\
 & = & \sum_{k=1}^{\infty}q^{k^{2}/2}(uq^{-\lambda})^{k}\\
 & + & \sum_{k=1}^{\infty}q^{k^{2}/2}(uq^{-\lambda})^{k}\left[(q^{m+k+1};q)_{\infty}-1\right]\\
 & = & \sum_{k=-\infty}^{-1}q^{k^{2}/2}(u^{-1}q^{\lambda})^{k}+r_{2}(n).\end{eqnarray*}
 By the binomial theorem, for $k\ge1$\begin{equation}
\left|(q^{m+k+1};q)_{\infty}-1\right|\le\frac{q^{m+2}(-q^{3};q)_{\infty}}{1-q},\label{eq:2.14}\end{equation}
 then\begin{equation}
|r_{2}(n)|\le\frac{q^{m+2}(-q^{3};q)_{\infty}\theta(\left|u\right|^{-1}q^{\lambda};q)}{1-q}.\label{eq:2.15}\end{equation}
 Thus we have proved that\begin{equation}
(-q^{-nt+1/2}u;q)_{\infty}=\frac{u^{m}\left\{ \theta(u^{-1}q^{\lambda};q)+r(n)\right\} }{(q;q)_{\infty}q^{m^{2}/2+m\lambda}}\label{eq:2.16}\end{equation}
 with\begin{eqnarray}
|r(n)| & \le & \frac{3(-q;q)_{\infty}\theta(\left|u\right|^{-1}q^{\lambda};q)}{1-q}\label{eq:2.17}\\
 & \times & \left\{ q^{m/2}+\frac{q^{m^{2}/8}}{|u|^{\left\lfloor m/2\right\rfloor +1}}\right\} .\nonumber \end{eqnarray}
 for $n$, $m$ and $\lambda$ satisfying (\ref{eq:2.2}) with $n$
and $m$ are large enough. 

In the case that $t$ is a positive irrational number, for any real
number $\beta\in[0,1)$, when $n$ and $m$ are large enough satisfying
(\ref{eq:2.5}) and (\ref{eq:2.6}), we have \begin{equation}
\beta+\gamma_{n}\ge-1,\label{eq:2.18}\end{equation}
 for these integers $n$, we take\begin{equation}
\nu_{n}=\left\lfloor -\frac{\log n}{\log q}\right\rfloor ,\label{eq:2.19}\end{equation}
 then we have \begin{eqnarray}
(-q^{-nt+1/2}u;q)_{\infty}(q;q)_{\infty} & = & \sum_{k=0}^{\infty}(q^{k+1};q)_{\infty}q^{k^{2}/2-km-k\beta-k\gamma_{n}}u^{k}\label{eq:2.20}\\
 & = & s_{1}+s_{2},\nonumber \end{eqnarray}
 with\begin{equation}
s_{1}=\sum_{k=0}^{m}(q^{k+1};q)_{\infty}q^{k^{2}/2-km-k\beta-k\gamma_{n}}u^{k}\label{eq:2.21}\end{equation}
 and\begin{equation}
s_{2}=\sum_{k=m+1}^{\infty}(q^{k+1};q)_{\infty}q^{k^{2}/2-km-k\beta-k\gamma_{n}}u^{k}.\label{eq:2.22}\end{equation}
 In $s_{1}$, we change summation from $k$ to $m-k$, \begin{eqnarray*}
\frac{s_{1}q^{mnt-m^{2}/2}}{u^{m}} & = & \sum_{k=0}^{m}(q^{m-k+1};q)_{\infty}q^{k^{2}/2}(u^{-1}q^{\beta+\gamma_{n}})^{k}\\
 & = & \sum_{k=0}^{\infty}q^{k^{2}/2}(u^{-1}q^{\beta})^{k}\\
 & - & \sum_{k=\nu_{n}+1}^{\infty}q^{k^{2}/2}(u^{-1}q^{\beta})^{k}\\
 & + & \sum_{k=0}^{\nu_{n}}q^{k^{2}/2}(u^{-1}q^{\beta})^{k}\left(q^{k\gamma_{n}}-1\right)\\
 & + & \sum_{k=0}^{\nu_{n}}q^{k^{2}/2}(u^{-1}q^{\beta+\gamma_{n}})^{k}\left\{ (q^{m-k+1};q)_{\infty}-1\right\} \\
 & + & \sum_{k=\nu_{n}+1}^{m}q^{k^{2}/2}(u^{-1}q^{\beta+\gamma_{n}})^{k}(q^{m-k+1};q)_{\infty}\\
 & = & \sum_{k=0}^{\infty}q^{k^{2}/2}(u^{-1}q^{\beta})^{k}+s_{11}+s_{12}+s_{13}+s_{14},\end{eqnarray*}
 thus there exists some constant $c_{11}(q,u)$ such that\begin{eqnarray*}
|s_{11}+s_{14}| & \le & 2\sum_{k=\nu_{n}+1}^{\infty}q^{k^{2}/2}\left|u\right|^{-k}\\
 & \le & c_{11}(q,u)\frac{\log n}{n}\end{eqnarray*}
 for $n$ large enough, and there exists some constant $c_{12}(q,u)$
such that\begin{eqnarray*}
|s_{12}+s_{13}| & \le & c_{12}(q,u)\frac{\log n}{n},\end{eqnarray*}
 for $n$ large enough. Thus for $n$, $m$ and $\beta$ satisfying
(\ref{eq:2.5}) and (\ref{eq:2.6}) and $n$ is large enough, there
exists some constant $c_{1}(q,u)$ such that\begin{equation}
\frac{s_{1}q^{mnt-m^{2}/2}}{u^{m}}=\sum_{k=0}^{\infty}q^{k^{2}/2}(u^{-1}q^{\beta})^{k}+e_{1}(n)\label{eq:2.23}\end{equation}
 with\begin{equation}
|e_{1}(n)|\le c_{1}(q,u)\frac{\log n}{n}.\label{eq:2.24}\end{equation}
 Similarly, in $s_{2}$ we change summation from $k$ to $k+m$, and
we could show that there for $n$, $m$ and $\beta$ satisfying (\ref{eq:2.5})
and (\ref{eq:2.6}) and $n$ is large enough, there exists some constant
$c_{2}(q,u)$ such that\begin{equation}
\frac{s_{2}q^{mnt-m^{2}/2}}{u^{m}}=\sum_{k=-\infty}^{-1}q^{k^{2}/2}(u^{-1}q^{\beta})^{k}+e_{2}(n)\label{eq:2.25}\end{equation}
 with\begin{equation}
|e_{2}(n)|\le c_{2}(q,u)\frac{\log n}{n}.\label{eq:2.26}\end{equation}
 Thus we have\begin{equation}
(-q^{-nt+1/2}u;q)_{\infty}=\frac{u^{m}\left\{ \theta(u^{-1}q^{\beta};q)+O\left(\frac{\log n}{n}\right)\right\} }{(q;q)_{\infty}q^{mnt-m^{2}/2}}\label{eq:2.27}\end{equation}
 for $n$, $m$ and $\beta$ satisfying (\ref{eq:2.5}) and (\ref{eq:2.6})
and $n$ is large enough. 
\end{proof}
\begin{conclusion*}
The phenomenon demonstrated here is universal for class of entire
basic hypergeometric functions, the general proof for this phenomenon
will be published elsewhere. In a separated work, we have proved the
same phenomenon happens with the Ramanujan's entire function $A_{q}(z)$,
or the $q$-Airy function.
\end{conclusion*}

\end{document}